\documentclass{amsart}
\usepackage{amsfonts,amssymb,amsmath,amsthm,latexsym}

\newcommand{\sC}{{\mathcal C}}

\newcommand{\sR}{{\mathcal R}}

\newcommand{\F}{{\mathbb{F}}}

\newtheorem{theorem}{Theorem}[section]
\theoremstyle{definition}
\newtheorem{definition}[theorem]{Definition}
\newtheorem{lemma}[theorem]{Lemma}

\newtheorem{proposition}[theorem]{Proposition}
\theoremstyle{remark}
\newtheorem{remark}[theorem]{Remark}

\newenvironment{defn*}{\begin{definition}}{\end{definition}}

\numberwithin{equation}{section}

\begin{document}

\title{NEW CUBIC SELF-DUAL CODES OF LENGTH 54, 60 AND 66}

\author{P{\i}nar {\c{C}omak}}
\address{Middle East Technical University, Department of Mathematics}
\email{pcomak@metu.edu.tr}

\author{Jon Lark Kim}
\address{Sogang University, Department of Mathematics}
\email{jlkim@sogang.ac.kr}

\author{Ferruh \"{O}zbudak}
\address{Middle East Technical University, Department of Mathematics and Institute of Applied Mathematics}
\email{ozbudak@metu.edu.tr}

\begin{abstract}
We study the construction of quasi-cyclic self-dual codes, especially of binary cubic ones. We consider the binary quasi-cyclic codes of length $3\ell$ with the algebraic approach of \cite{LS01}. In particular, we improve the previous results by constructing 1 new binary $[54,27,10]$, 6 new $[60,30,12]$ and 50 new $[66,33,12]$ cubic self-dual codes. We conjecture that there exist no more binary cubic self-dual codes with length 54, 60 and 66.
\end{abstract}

\maketitle

\section{Introduction}
Quasi-cyclic and self-dual codes are interesting classes of linear codes. Quasi-cyclic codes are those linear codes which takes the maximum possible value of minimum distance among the codes with the same length and same dimension. Another class of interesting linear codes is the classes of self-dual codes. Self-dual codes have close connections with group theory, lattice theory and design theory. There has been an active research on the classification of self-dual codes over finite fields and over rings. 
In \cite{LS01} it was shown that all cubic binary codes of length $3\ell$ can be obtained by a generalization of cubic construction of Turyn's, from a binary code and a quaternary code of both length $\ell$.

The rest of this paper is organized as follows. In Section 2, we give some preliminaries about linear codes. In Sections 3 and 4, we recall the basic properties of quasi-cyclic codes and the cubic construction that is used. In Section 5, we present our results.

\section{Preliminaries}
A $q$-ary linear code $\sC$ is a linear subspace of $\F_q^n$. If $\sC$ has dimension $k$, then $\sC$ is called an $[n,k]$-linear code or $[n,k,d]$-linear code, where $d=d(\sC)$ is the minimum Hamming distance which is the minimum number of distinct coordinates between any pair of distinct codewords in $\sC$. The Hamming weight $w(c)$ of a codeword $c$ in $\sC$ is defined to be the number of non-zero entries of $c$. For a linear code, we have that $d(\sC)=w(\sC)$. Two codes are said to be equivalent up to permutation if they differ only in the order of their coordinates. The Hamming weight enumerator of the code $\sC$ is defined to be $ W_{\sC}(y)= \sum_{c\in \sC}{y^{wt(c)}}=\sum_{i=0}^{n}{A_i y^{i}}$, where $A_i$ is the number of vectors of the code $\sC$ having Hamming weight $i$. 

A linear code $\sC$ is said to be cyclic if for every codeword $c = (c_0,c_1,.., c_{n-1}) \in \sC$ there exists a corresponding codeword $c'= (c_{n-1},c_0,\dotsc, c_{n-2}) \in \sC$, in which case we call $c'$ is a cyclic shift of $c$. It is more convenient to represent the codeword $c = (c_0,c_1,\dotsc, c_{n-1})$ as the polynomial $c(x)=c_0 + c_1x +c_2x^2+\dotsb+c_{n-1}x^{n-1}$, where $c_i\in \F_q$, $\forall i={0,..,n-1}$. In this representation, the cyclic shift of $c(x)$ is simply $c'(x) = xc(x)$. Since any linear combination of cyclic shifts of $c(x)$ is also a codeword of $\sC$, a cyclic code is an ideal in $\F_q[x]/(x^n-1)$. In a nonzero cyclic code $\sC$, the unique monic polynomial of smallest degree is called the generator polynomial of $\sC$. The generator polynomial of a cyclic $[n,k]$-code has degree $n-k$ and denoted by $g(x)=g_0+g_1x+\dotsb+g_{n-k}x^{n-k}$, which is a factor of $x^n-1$. The cyclic code $\sC$ is shown as $\sC=\langle g(x) \rangle$. For a cyclic $[n,k]$-code $\sC$, a matrix $G \in \F_q^{k \times n}$ is a generator matrix for $\sC$ if its $k$ rows span $\sC$. 

We define the dual of a code $\sC$ to be $\sC^{\perp} = \{u \in \F_q^n : (u,v)=0 \text{ for all } v\in \sC\}.$ Here the inner product is the standard (Euclidean) inner product. The dual code $\sC^{\perp}$ of the code $\sC$ over $\F_q$ and over $\sR:=\sR(\F_q,m)=\F_q[Y]/(Y^m-1)$ is understood with respect to the Euclidean (standard) inner product and the Hermitian inner product, respectively. If $\sC$ is an $[n,k]$-code over $\F_q$, then the dual code $\sC^{\perp}$ is a linear $[n,n-k]$ code. For any $[n,k]$-code $\sC$, we have $(\sC^{\perp})^{\perp}=\sC$. A code $\sC$ is said to be self-orthogonal if $\sC \subset{\sC^{\perp}}$ and self-dual if $\sC=\sC^{\perp}$. If a code $\sC$ of length $n$ is self-dual, then $n$ must be even; and $\sC$ is a subspace of dimension $n/2$. 

If $\sC \subset \F_2^n$ is a binary self-dual code, then the weight of all codewords must be even. The binary self-dual codes in which there is at least one codeword with weight not divisible by 4 are called Type I or singly-even self-dual binary codes. Otherwise, the binary self-dual codes are called Type II or doubly-even self-dual binary codes.

\subsection{Inner Products}
In order to define dual codes, we need to define inner products. Inner products for linear codes over $\F_q^{\ell m}$ and over $\sR^{\ell}$ are defined as follows \cite{HA98}.

\subsubsection{Euclidean inner product} 
Defined on $\F_q^{\ell m}$ as $$ a \cdot b = \sum_{i=0}^{m-1}\sum_{j=0}^{\ell-1}a_{ij}b_{ij} $$ for $a=(a_{0,0},a_{0,1},\dotsc,a_{0,\ell-1},a_{1,0},\dotsc,a_{1,\ell-1},\dotsc,a_{m-1,0},\dotsc,a_{m-1,\ell-1})$ and \\ $b=(b_{0,0},b_{0,1},\dotsc,b_{0,\ell-1},b_{1,0},\dotsc,b_{1,\ell-1},\dotsc,b_{m-1,0},\dotsc,b_{m-1,\ell-1})$. 

\subsubsection{Hermitian inner product} 
Defined on $\sR^\ell$ as $$ \langle x,y \rangle=\sum_{j=0}^{\ell-1}x_j\overline{y_j}$$ for $ x=(x_0,x_1,\dotsc,x_{\ell-1}) $ and $ y=(y_0,y_1,\dotsc,y_{\ell-1})$. Here the conjugation operation ${}^{-}$ on $\sR$ sends $Y$ to $Y^{-1}=Y^{m-1}$ (identity on $\F_q$) and satisfies $\overline{\overline{x}}=x, \overline{x+y} = \overline{x} + \overline{y}, \overline{xy} = \overline{x}\text{ }\overline{y}.$

\section{Quasi-Cyclic Codes}
Let $\F_q$ be a finite field and $m$ be a positive integer coprime with the characteristic of $\F_q$. A linear code $\sC$ of length $\ell m$ over $\F_q$ is called quasi-cyclic code if $(c_{m-1,0},\dotsc,c_{m-1,\ell-1},c_{0,0},\dotsc,c_{0,\ell-1},\dotsc,c_{m-2,0},\dotsc,c_{m-2,\ell-1}) \in \sC$ whenever the codeword $(c_{0,0},\dotsc,c_{0,\ell-1},c_{1,0},\dotsc,c_{1,\ell-1},\dotsc,c_{m-1,0},\dotsc,c_{m-1,\ell-1}) \in \sC $. 

This code is invariant under $\ell$-shift and such codes are called as $\ell$-quasi-cyclic codes or quasi-cyclic codes of index $\ell$. The quasi-cyclic codes are the generalization of cyclic codes. Cyclic codes correspond to the case $\ell=1$. 

\subsection{1-1 correspondence}
Let $\F_q[Y]$ denote the polynomial ring over $\F_q$. Consider the ring $\sR:=\sR(\F_q,m)=\F_q[Y]/(Y^m-1)$. 
This ring is the same as the polynomial representation of cyclic codes of length $m$ over $\F_q$. Namely, the cyclic codes of length $m$ over $\F_q$ are ideals of $\sR(\F_q,m)$. Modules over $\sR$ is closely related to the ideals in $\F_q[Y]/(Y^m-1)$. In fact, ideals are just 1-dimensional $\sR$-submodules. Let $\sC$ be a $\ell$-quasi-cyclic code over $\F_q$ of length $\ell m$ and let $$c=(c_{0,0},\dotsc,c_{0,\ell-1},c_{1,0},\dotsc,c_{1,\ell-1},\dotsc,c_{m-1,0},\dotsc,c_{m-1,\ell-1})$$ denote a codeword in $\sC$. Define a map $\phi:\F_q^{\ell m} \rightarrow \sR^\ell$ by $$
\phi(c)=(c_0(Y),c_1(Y),\dotsc,c_{\ell-1}(Y)) \in \sR^{\ell}
$$
where 
$c_j(Y)=\sum_{i=0}^{m-1}c_{ij}Y^i \in \sR, \quad j=0,\dotsc,\ell-1$. 

A linear code $\sC$ of length $n$ over $\sR$ is defined to be a $\sR$-submodule of $\sR^n$.
If $\sR$ is a finite field $\F_q$ of order $q$, the linear code $\sC$ of order $n$ over $\F_q$ is an $\F_q$-vector subspace of $\F_q^n$.

\begin{lemma} {\cite{LS01}}
The map $\phi$ gives a one-to-one correspondence between $\ell$-quasi-cyclic codes over $\F_q$ of length $\ell m$ and linear codes over $\sR$ of length $\ell$.
\end{lemma}

\begin{proposition} \label{pr:2.2} \cite{LS01}
Let $a,b \in \F_q^{\ell m}$. Then
$(\mathit{T}^{\ell k}(a))\cdot b =0$ for $0\leq k \leq m-1$ if and only if $\langle\phi(a),\phi(b)\rangle=0$.
\end{proposition}

It follows from Proposition \ref{pr:2.2} that a quasi-cyclic code $\sC$ is self-dual with respect to the Euclidean inner product if and only if $\phi(\sC)$ is self-dual with respect to the Hermitian inner product, where $\sC$ is an $\ell$-quasi-cyclic code over $\F_q$ of length $\ell m$ and $\phi(\sC)$ is its image in $\sR^\ell$ under $\phi$. We also have that ${\phi(\sC)}^\perp={\phi(\sC^\perp)}$, where the dual in $\F_q^{\ell m}$ is taken with respect to the Euclidean inner product and the dual in $\sR^\ell$ is taken with respect to the Hermitian inner product.

\subsection{Existence of Self-Dual Codes}
In \cite{JLK12}, it is proved that there exist self-dual binary codes of length $\ell$ over $\sR=\sR(\F_2,m)=\F_2[Y]/(Y^m-1)$ if and only if $2 | \ell$. For binary $\ell$-quasi-cyclic self-dual codes of length $\ell m$, if $m$ is a prime not dividing $i$, then $m$ must divide $A_i$, which is the number of codeword with Hamming weight $i$. This gives the possible weight enumerators of self-dual codes of a given length.

\subsection{Ring Decomposition}
Let $\sR=\sR(\F_q,m)=\F_q[Y]/(Y^m-1).$ If $\gcd(m,q)=1$, then the ring can be decomposed into a direct sum of fields by the Chinese remainder theorem (CRT) or discrete Fourier transform (DFT) \cite{LS01}. By this approach, the quasi-cyclic codes can be decomposed into codes of lower lengths. The polynomial $Y^m-1$ factors completely into distinct irreducible factors in $\F_q[Y]$ as
\begin{equation} \label{eq:factor}
Y^m-1=\delta g_1 \dotsc g_s h_1 h_1^* \dotsc h_t h_t^*
\end{equation}
where $\delta$ is nonzero in $\F_q$, $g_1 \dotsc g_s$ are the polynomials which are self-reciprocal, and $h_i^*$'s are reciprocals of $h_i$'s, for all $1 \leq i \leq t$. Then by CRT \cite{LS01}, the ring $\sR$ can be written as
\begin{equation} \label{eq:decomposition}
\sR=\frac{F_q[Y]}{(Y^m-1)}=\Bigg( \overset{s}{\underset{i=1}{\bigoplus}} \frac{F_q[Y]}{(g_i)} \Bigg) \oplus
 \Bigg( \overset{t}{\underset{j=1}{\bigoplus}} \bigg( \frac{F_q[Y]}{(h_j)} \oplus \frac{F_q[Y]}{(h_j^*)} \bigg) \Bigg).
\end{equation}

For notational convenience, let $G_i$, $H_j'$ and $H_j''$ denote $F_q[Y] \big/(g_i)$, $F_q[Y] \big/ (h_j)$ and $F_q[Y] \big/(h_j^*)$, respectively. Then 

$$ 
\sR^\ell = \Bigg( \overset{s}{\underset{i=1}{\bigoplus}} G_i^\ell \Bigg) \oplus \Bigg( \overset{t}{\underset{j=1}{\bigoplus}} \bigg( H_j'^\ell \oplus H_j''^\ell \bigg) \Bigg).
$$

and every $\sR$-linear code $\sC$ of length $\ell$ can be decomposed as the direct sum

$$
\sC = \Bigg( \overset{s}{\underset{i=1}{\bigoplus}} \sC_i \Bigg) \oplus \Bigg( \overset{t}{\underset{j=1}{\bigoplus}} \bigg( \sC_j' \oplus \sC_j'' \bigg) \Bigg)
$$
where $\sC_i$, $\sC_j'$ and $\sC_j''$ are linear codes over $G_i$, $H_j'$ and $H_j''$, respectively, all of length $\ell$ for each $1 \leq i \leq s$, and for each $1 \leq j \leq t$. 

Let $x=(x_0,x_1,\dots,x_{\ell-1})$ and $y=(y_0,y_1,\dots,y_{\ell-1})$. Here, for $1 \leq i \leq s$, the Hermitian inner product of $x$ and $y$ with $x_i, y_i \in G_i$ is defined as in \cite[Section IV]{LS01} by $\langle x,y \rangle = x_0y_0^{m-1}+ \cdots +x_{\ell-1}y_{\ell-1}^{m-1} $. Moreover, for $1 \leq i \leq t$, the Euclidean inner product of $x$ and $y$ with $x_i,y_i\in H_j'$ is defined by $x \cdot y = x_0 y_0 + \cdots + x_{\ell-1}y_{\ell-1}$.

Notice that every element of $\sR$ can be written as $c(Y)$, for some polynomial $c\in \F_q[Y]$. The decomposition 
(\ref{eq:decomposition}) shows that $c(Y)$ may also be written as an $(s+2t)$-tuple
\begin{equation} \label{eq:decomp}
(c_1(Y),\dotsc,c_s(Y),c_1'(Y),c_1''(Y),\dotsc,c_t'(Y),c_t''(Y))
\end{equation}
where $c_i(Y) \in G_i$, $(1 \leq i \leq s)$, $c_j'(Y) \in H_j'$, $(1 \leq j \leq t)$ and $c_j''(Y) \in H_j''$, $(1 \leq j \leq t)$. Here, we can consider the $c_i,c_j'$ and $c_j''$ as polynomials in $\F_q[Y]$. We can also write 
$$\overline{c(Y)}=(\overline{c_1(Y)},\dotsc,\overline{c_s(Y)},c_1''(Y),c_1'(Y),\dotsc,c_t''(Y),c_t'(Y)).$$

\begin{theorem} \label{thm:decomperp}  \emph{\cite{LS01}} 
An $\ell$-quasi-cyclic code $\sC$ of length $\ell m$ over $\F_q$ is self-dual if and only if 
$$
\sC = \Bigg( \overset{s}{\underset{i=1}{\bigoplus}} \sC_i \Bigg) \oplus \Bigg( \overset{t}{\underset{j=1}{\bigoplus}} \bigg( \sC_j' \oplus (\sC_j')^{\perp} \bigg) \Bigg)
$$
where, for $1 \leq i \leq s$, $\sC_i$ is a self-dual code over $G_i$ of length $\ell$ with respect to the Hermitian inner product and for $1 \leq j \leq t$, $\sC_j'$ is a linear code of length $\ell$ over $H_j'$ and $(\sC')^{\perp}$ is its dual with respect to the Euclidean inner product as defined above.
\end{theorem}

\section{Cubic Self-Dual Binary Codes and Cubic Construction}
There are some construction methods for combining codes to get new codes with greater length for different values of $q$, $m$ and $\ell$ (for example see \cite{BO03}).

In this work, we focus on the case of so called binary cubic codes, which ,s the case for $q=2$ and $m=3$. We use a cubic construction in \cite{BO03} and \cite{LS01} to find new codes. 

Since $Y^2 + Y + 1$ is irreducible in $\F_2[Y]$, we can write $Y^3 - 1 = (Y-1)(Y^2 + Y + 1)$ as a product of irreducible factors. By (\ref{eq:decomposition}), $\sR$ can be decomposed as
$$
\sR = \frac{\F_2 [Y]}{(Y^3 -1)} = \F_2 \oplus \F_{2^2}.
$$
This gives a correspondence between the $\ell$-quasi-cyclic codes $\sC$ of length $3 \ell$ over $\F_2$ and a pair $(\sC_1, \sC_2)$, where $\sC_1$ is a linear code over $\F_2$ of length $\ell$ and $\sC_2$ is a linear code over $\F_4 = \F_2(\omega)$ of length $\ell$ where $\omega^2 + \omega + 1 = 0$. Using the discrete Fourier transform \cite{LS01}, we can define a Gray map from $\F_2^{\ell} \times \F_4^{\ell}$ to $\F_2^{3\ell}$ by
\begin{equation} \label{eq:cubic}
\sC = \phi(\sC_1,\sC_2) = \{(x + a\  |\  x + b\  |\  x + a + b)\ |\  x \in \sC_1,\ a + \omega b \in \sC_2 \}
\end{equation}

Moreover, $\sC$ is self-dual if and only if $\sC_1$ is self-dual with respect to the Euclidean inner product and $\sC_2$ is self-dual with respect to the Hermitian inner product.

It is shown in \cite{LS01} that all such codes can be obtained by this method, from a binary code over $\F_2$ and a quaternary code over $\F_4$ both of length $\ell$. Here, as in \cite{BO03}, cubic binary codes of length $3\ell$ are viewed as codes of length $\ell$ over the ring $\F_2 \times \F_{2^2}$. 

It is easily verified that, if $a,b \in \sC_2'$ for some linear code $\sC_2'$ over $\F_2$, then $\sC_2 = \{a + b \omega \ | \ a,b \in \sC_2' \}$ is a linear code over $\F_{2^2}$, where $\omega^2+\omega+1=0$.
So, if we begin with two $\F_2$-linear codes $\sC_2'$ and $\sC_1$, the construction in (\ref{eq:cubic}) gives Turyn's $(a + x\ | \ b + x \ | \ a + b + x)$-construction. In particular, we obtain the following.

\begin{theorem} \emph{\cite{LS01}} \label{thm:constr}
The $(x + a \ | \ x + b \ | \ x + a + b)$-construction which is applied to two linear codes $\sC_1$ and $\sC_2'$ over $\F_{2}$ of length $\ell$, gives an $\F_{2}$-linear code $\sC$ of length $3 \ell$ that is quasi-cyclic of index $\ell$.
\end{theorem}

The authors of \cite{BO03}, \cite{BO06}, \cite{JLKBO12} and \cite{JLK12} completed the classification of binary cubic self-dual codes of lengths $\leq 48$ up to permutation equivalence. The number of cubic self-dual codes in various cases is given as follows:
\begin{itemize} 
	\item[(i)] for $\ell=2$, unique binary cubic self-dual code of length 6, 
	\item[(ii)] for $\ell=4$, 2 binary cubic self-dual codes of length 12, 
	\item[(iii)] for $\ell=6$, 3 binary cubic self-dual codes of length 18,
	\item[(iv)] for $\ell=8$, 16 binary cubic self-dual codes of length 24,
	\item[(v)] for $\ell=10$, 8 binary cubic self-dual codes of length 30, 
	\item[(vi)] for $\ell=12$, 13 binary cubic self-dual codes of length 36,
	\item[(vii)] for $\ell=14$, 1569 binary cubic self-dual codes of length 42,
	\item[(viii)] for $\ell=16$, 264 binary cubic self-dual codes of length 48,
	\item[(ix)] for $\ell=18$, $\geq 13$ binary cubic self-dual codes of length 54,
	\item[(x)] for $\ell=20$, $\geq 3$ binary cubic self-dual codes of length 60,
	\item[(xi)] for $\ell=22$, $\geq 7$ binary cubic self-dual codes of length 66.
\end{itemize}

\section{Main Results}
\subsection{$\ell=18$, $[54,27,10]$ codes:}
The shortest length of binary cubic self-dual codes for which the classification is not completed is 54. The number of known inequivalent codes that were found in \cite{BO06} is 13. In this paper, we find 1 more such code by the cubic construction (\ref{eq:cubic}). 

For the self-dual $[54,27,10]$ codes, the following are the only possible weight enumerators \cite{CS90}:
\begin{small}
\begin{flalign} \nonumber 
W_1 &= 1 + (351 - 8 \beta) y^{10} + (5031 + 24 \beta) y^{12} + (48492 + 32 \beta) y^{14} + \ldots \quad 0 \leq \beta \leq 43 \\\nonumber
W_2 &= 1 + (351 - 8 \beta) y^{10} + (5543 + 24 \beta) y^{12} + (43884 + 32 \beta) y^{14} + \ldots \quad 12 \leq \beta \leq 43.
\end{flalign}
\end{small}
\begin{lemma} \label{la:div} \cite{JLK12}
Let $\sC$ be a binary $\ell$-quasi-cyclic self-dual code of length $m \ell$ where $m$ is prime. If $m$ does not divide the weight $i$, then $m$ must divide $A_i$, where $ W_{\sC}(y)= \sum_{c\in \sC}{y^{wt(c)}}=\sum_{i=0}^{n}{A_i y^{i}}$.
\end{lemma}

In \cite{BO06}, seven inequivalent codes with $W_1$ for $\beta = 0,3,6,9,12,15,18$ and six inequivalent codes with $W_2$ for $\beta = 12,15,18,21,24,27$ are found. Note that $\beta$ should be divisible by 3 by the above lemma. 

By the construction (\ref{eq:cubic}), binary codes $\sC$ of length 54 are formed from a binary code $\sC_1$ of length 18 and a quaternary code $\sC_2$ of length 18. In \cite{BO03}, the following proposition is proved.

\begin{proposition} \label{pro:dis} \cite{BO03}
If $\sC = \phi(\sC_1, \sC_2)$, then $d(\sC) \leq \min(3d(\sC_1),2d(\sC_2))$.
\end{proposition}

It is known that there are 7 binary $[18,9,2]$, 2 binary $[18,9,4]$ codes and 244 quaternary $[18,9,6]$ and 1 quaternary $[18,9,8]$ codes which are all self-dual and listed in \cite{MH16}. To obtain more $[54,27,10]$ codes, by the help of proposition \ref{pro:dis}, it is sufficient to take the only 2 binary codes with $d=4$ and all quaternary self-dual codes with a huge number of permutations.

By the construction above, we found eight $[54,27,10]$ codes with weight enumerator $W_1$ for $\beta = 0,3,6,9,12,15,18,21$ and six $[54,27,10]$ codes with weight enumerator $W_2$ for $\beta = 12,15,18,21,24,27$ by taking the 2 binary codes for $\sC_1$ and $18^{th}$ and $38^{th}$ $[18,9,6]$ self-dual quaternary codes taken from \cite{MH16} for $\sC_2$.

These codes are of Type I 18-quasi-cyclic self-dual codes of length 54 since their binary components are of Type I and self-dual with respect to the Euclidean inner product.

\begin{remark} \label{re:conj1} Based on computational evidence, we conjecture that there are no other $[54,27,10]$ self-dual cubic codes over $\F_2$.
\end{remark}

Our computational results, with $\beta$ a multiple of 3, are listed below: 
\renewcommand{\arraystretch}{1.5}
\begin{flushleft}
\begin{tabular*}{0.75\textwidth}{r|c|c|c|}
\multicolumn{1}{r}{}
 & \multicolumn{1}{c}{Possible values}
 & \multicolumn{1}{c}{Values Thm.\ref{eq:cubic}} 
 & \multicolumn{1}{c}{Conjecture, Rk.\ref{re:conj1}} \\
\cline{2-4}
$W_1$ & $0 \leq \beta \leq 43$ & $\beta \in \{ 0,3,6,9,12,15,18,\textbf{21}\}$ & $\beta \notin \{24,27 \cdots,42 \}$ \\
\cline{2-4}
$W_2$ & $12 \leq \beta \leq 43$ & $\beta \in \{ 12,15,18,21,24,27 \}$ & $\beta \notin \{30,33 \cdots,42 \}$ \\
\cline{2-4}
\end{tabular*}
\end{flushleft}

\subsection{$\ell =20$, $[60,30,12]$ codes:}
In this case, the number of previously known inequivalent codes was 3. In this paper, we find 6 more such codes by the cubic construction (\ref{eq:cubic}). 

For self-dual $[60,30,12]$ codes, the following are the only possible weight enumerators \cite{GH96}: \begin{flalign}
 \nonumber W_1 &= 1 + 2555 y^{12} + 33600 y^{14} + 278865 y^{16}+ \ldots &&\\\nonumber
 W_2 &= 1 + 2619 y^{12} + 33216 y^{14} + 279441 y^{16}+ \ldots \\\nonumber
 W_3 &= 1 + 3195 y^{12} + 29760 y^{14} + 284625 y^{16}+ \ldots \\\nonumber
 W_4 &= 1 + 3451 y^{12} + 24128 y^{14} + 336081 y^{16}+ \ldots
\end{flalign}

By Lemma \ref{la:div}, we know that there is no code with weight enumerator $W_4$ since $3 \nmid 14$ and $3 \nmid 24128$. 

In \cite{BO03}, three inequivalent extremal codes with $W_3$ are found with automorphism groups of size 3, 6 and 12. 

By the construction (\ref{eq:cubic}), binary codes $\sC$ of length 60 are formed from a binary code $\sC_1$ of length 20 and a quaternary code $\sC_2$ of length 20. 

It is known that there are 9 binary self-dual $[20,10,2]$ codes, 7 binary self-dual $[20,10,4]$ codes and 245 quaternary self-dual $[20,10,2]$ codes, 2181 quaternary self-dual $[20,10,4]$ codes, 999 quaternary self-dual $[20,10,6]$ codes and 2 quaternary self-dual $[20,10,8]$ codes listed in \cite{MH16}. To obtain $[60,30,12]$ codes, since $d(\sC) \leq \min(3d(\sC_1),2d(\sC_2))$, it is sufficient to take only the binary codes of minimum distance 4 and the quaternary codes of minimum distance 6 and 8 with a huge number of permutation.

By this construction, we found nine inequivalent $[60,30,12]$ codes with weight enumerator $W_3$ with automorphism groups of size 3,6,9,12,18,24,30,48,60. 

These $[60,30,12]$ codes are of Type I 20-quasi-cyclic self-dual codes since their binary components $\sC_1$'s are of Type I and self-dual with respect to the Euclidean inner product.

\begin{remark} \label{re:conj2} Based on computational evidence, we conjecture that there are no other $[60,30,12]$ self-dual cubic codes over $\F_2$.
\end{remark}

Our computational results are listed below: 
\renewcommand{\arraystretch}{1.5}
\begin{flushleft}
\begin{tabular*}{0.75\textwidth}{ r|c|c|c| }
\multicolumn{1}{r}{}
 & \multicolumn{1}{c}{Found values} 
 & \multicolumn{1}{c}{New values} 
 & \multicolumn{1}{c}{Conjecture, Rk.\ref{re:conj2}} \\
\cline{2-4}
$W_3$ & $3,6,9,12,18,24,30,48,60$ & $ \textbf{9,18,24,30,48,60}$ & no other codes \\
\cline{2-4}
\end{tabular*}
\end{flushleft}

\subsection{$\ell =22$, $[66,33,12]$ codes:}
The number of previously known inequivalent codes was 7 \cite{JLK12}. In this paper, we find 50 more such codes by the cubic construction (\ref{eq:cubic}). 

For self-dual $[66,33,12]$ codes, there are three possible weight enumerators:
\begin{flalign}
\nonumber W_1 &= 1 + (858 + 8 \alpha ) y^{12} + (18678 + 24 \alpha) y^{14} + \ldots \quad 0 \leq \alpha \leq 778 &&\\\nonumber
 W_2 &= 1 + (858+ 8 \alpha) y^{12} + (18166 + 24\alpha) y^{14} + \ldots \quad 14 \leq \alpha \leq 756 \\\nonumber
 W_3 &= 1 + 1690 y^{12} + 7990 y^{14} + \ldots
\end{flalign}

Before our work, in \cite{JLK12} it was known that seven inequivalent codes with $W_1$ for $\alpha = 17, 21, 23, 26, 30, 43, 46$ are found.

By Lemma \ref{la:div}, there is no code with weight enumerators $W_2$ and $W_3$ since $3 \nmid (18166 + 24\alpha)$ and $3 \nmid 7990$. Therefore any binary cubic self-dual [66, 33, 12] code should have weight enumerator $W_1$.

By the construction (\ref{eq:cubic}), binary codes $\sC$ of length 66 are formed from a binary code $\sC_1$ of length 22 and a quaternary code $\sC_2$ of length 22. 

It is known that there are 16 binary self-dual $[22,11,2]$ codes, 8 binary self-dual $[22,11,4]$ codes, 1 binary self-dual $[22,11,6]$ code and 723 quaternary self-dual $[22,11,8]$ codes listed in \cite{MH16}. To obtain $[66,33,12]$ codes, by proposition \ref{pro:dis}, we take the binary codes of minimum distance 4 and 6, and the quaternary codes of minimum distance 8 with a huge number of permutation.

By this construction, we found fifty seven inequivalent $[66,33,12]$ codes with weight enumerator $W_1$ for $\alpha = 6,8-54,56,57,59,60,62,65,68,69,71$. 

These codes are of Type I 22-quasi-cyclic self-dual codes of length 66 since their binary components $\sC_1$'s are of Type I and self-dual with respect to the Euclidean inner product.

\begin{remark} \label{re:conj3} Based on computational evidence, we conjecture that there is no other $[66,33,12]$ self-dual cubic code over $\F_2$.
\end{remark}

Our computational results are listed below: 
\renewcommand{\arraystretch}{1.5}
\begin{flushleft}
\begin{tabular*}{0.95\textwidth}{ r|c|c|c| }
\multicolumn{1}{r}{}
 & \multicolumn{1}{c}{Possible values}
 & \multicolumn{1}{c}{Previous values} 
 & \multicolumn{1}{c}{Found values} \\
\cline{2-4}
$W_1$ & $0 \leq \alpha \leq 778$ & 17,21,23,26,30,43,46 & \begin{tabular}{lp{.2\textwidth}l} 6,8-54,56,57,59,60,\\62,65,68,69,71 \end{tabular}\\
\cline{2-4}
\end{tabular*}
\end{flushleft}

\bibliographystyle{amsplain}
\nocite{*}

\bibliography{references}
 
\end{document}